\documentstyle[12pt,amssymb]{amsart}
\newcommand{\stopproof}{\hfill \nobreak\medskip $\blacksquare$ \\
\hspace*{\fill}}
\newcommand{\dom}{\mbox{\rm dom}}
\newcommand{\cof}{\mbox{\rm cof}}

\newcommand{\AND}{\mbox{ \rm and }}

\newcommand{\IF}{\mbox{ \rm if }}

\newcommand{\forces}[2]{\Vdash_{#1} \mbox{``} #2 \mbox{''}}

\newcommand{\proof}{{\bf Proof:} \ }

\newcommand{\PP}{{\Bbb P}}
\newcommand{\alpeh}{{\aleph}}

\newcommand{\DD}{{\Bbb D}}
\newcommand{\RR}{{\Bbb R}}
\newcommand{\CC}{{\Bbb C}}

\newcommand{\QQ}{{\Bbb Q}}
\newcommand{\SS}{{\Bbb S}}

\newcommand{\presup}[2]{\, ^{#1} \! #2}
\newcommand{\fomom}{\presup{\omega}{\omega}}
\newcommand{\wfomom}{\presup{\stackrel{\omega}{\smile}}{\omega}}
\newcommand{\restricts}{\! \upharpoonright \!}
\newcommand{\lcard}{\, \mid \!}
\newcommand{\rcard}{\! \mid \,}

\newcommand{\card}[1]{\lcard #1 \rcard}
\newtheorem{theor}{Theorem}[section]

\newtheorem{defin}{Definition}[section]
\newtheorem{corol}{Corollary}[section]
\newtheorem{lemma}{Lemma}[section]

\title{Maximal Chains in $\fomom$ and Ultrapowers
of the Integers}
\author{Saharon Shelah}
\address{Institute of Mathematics \\ Hebrew University \\ Jerusalem,
Givat Ram,  Israel and Department of Mathematics\\
Rutgers University\\ New Brunswick, New Jersey}
\author{Juris Stepr\={a}ns}
\address{Department of Mathematics, York University \\
4700 Keele Street \\ North York, Ontario \\ Canada \ \ \ \ M3J 1P3}                              
\thanks{The first author is partially  supported by the basic reasearch
fund of the Israeli Academy. The second author is partially
supported by NSERC and was a guest of Rutgers University while the
research on this paper was being done. The authors would also like to thank
P. Nyikos for his valuable comments on early versions of this paper.
This is number 465 on the first author's list of publications}
\begin{document}
\maketitle
\bibliographystyle{amsplain}
\begin{abstract}Various questions posed by P. Nyikos 
concerning ultrafilters on $\omega$ and chains in the partial order
$(\omega,<^*)$ are answered. The main tool is the oracle chain condition
and variations of it.

\end{abstract}
\section{Introduction}
In \cite{nyik.c} various axioms related to maximal chains in ultrapowers
 of the integers were classified and studied. The purpose of
the present paper is to answer several of the questions posed in
that paper and to pose some new ones. 

The notation and terminology of this 
paper will adhere as much as possible to accepted
standards but some of the main points are listed here.
 The relation $a \subset^{\ast} b$ means that $\card{a
\setminus b}  < \aleph_0$ while $f\leq^* g$ means that $f$ and $g$ belong
to $\fomom$ --- or, perhaps, $\presup{A}{\omega}$ where
$A\subseteq\omega$ is infinite --- and $f(n) \leq g(n)$ for all but
finitely many integers $n$. If $f(n) < g(n)$ for all but
finitely many integers $n$ then this will be denoted by $f<^*g$. 
 By a {\em chain} in $\fomom$ will be meant a
subset of $\fomom$ which is well ordered by $<^*$ and consists of
nondecreasing functions. In the next section the effects of modifying
 this definition of a chain
will be discussed. A subset
$\cal{S} \subseteq \fomom$ will be said to be unbounded if  for every
 $f\in\fomom$ there is $g \in \cal{S}$ such that $g\not\leq^* f$.
The least cardinality of an unbounded subset of $\fomom$ is denoted
by $\frak{b}$ while 
the least cardinality of a cofinal subset of $\fomom$ is denoted
by $\frak{d}$.
 The term {\em ultrafilter} will be reserved
for   ultrafilters on $\omega$ which contain no finite sets. 
A
$P$-point  is an ultrafilter on $\omega$, $\cal{U}$,
such that for every ${\cal A} \in [\cal{U}]^{\aleph_0}$ there is $B \in 
\cal{U}$
such that $B \subset^{\ast} A$ for every $A \in \cal A$. If
$\cal{U}$ is a filter then $\cal{U}^{\ast}$ will denote the dual
ideal to $\cal{U}$.  

If $\cal{U}$  is an ultrafilter then
{\em the integers modulo $\cal{U}$} will refer to the ultrapower of
the integers with respect to $\cal{U}$ and will be denoted by
$\fomom/\cal{U}$.  If $\cal{U}$ is an ultrafilter then 
 $\cal{C}\subseteq\fomom$ will be said to be unbounded modulo
$\cal{U}$ if, letting $[f]_{\cal{U}}$ represent the equivalence
class of $f\in\fomom$ in $\fomom/\cal{U}$, the set $\{
[f]_{\cal{U}} : f \in \cal{C}\}$ is unbounded in the linear order
$\fomom/\cal{U}$. The least ordinal which can be embedded 
cofinally in a
linear ordering  $L$  
is denoted by $\cof(L)$ --- $\cof(\fomom/\cal{U})$  will be an
important invariant of $\cal{U}$ in the following discussion.

For reference, here are 
Nyikos' axioms (throughout $\cal{C}$ refers to a maximal chain
of nondescending functions in $\fomom$ and $\cal{U}$ refers to an ultrafilter)
\begin{itemize}
\item {\bf Axiom~1\ \ \ } $(\forall \cal{U})(\forall \cal{C})(\cal{C}
\mbox{ is unbounded modulo }\cal{U})$
    \item {\bf Axiom~2\ \ \ } $(\exists \cal{C})(\forall \cal{U})(\cal{C}
\mbox{ is unbounded modulo }\cal{U})$
\item {\bf Axiom~3\ \ \ } $(\exists \cal{U})(\forall \cal{C})(\cal{C}
\mbox{ is unbounded modulo }\cal{U})$
\item {\bf Axiom~4\ \ \ } $(\forall \cal{U})(\exists \cal{C})(\cal{C}
\mbox{ is unbounded modulo }\cal{U})$
\item {\bf Axiom~5\ \ \ } $(\forall \cal{C})(\exists \cal{U})(\cal{C}
\mbox{ is unbounded modulo }\cal{U})$
\item {\bf Axiom~5.5\ } $(\exists \cal{U})(\cof(\fomom/\cal{U}) = \frak{b})$
\item {\bf Axiom~6\ \ \ } $(\exists \cal{C})(\exists \cal{U})(\cal{C}
\mbox{ is unbounded modulo }\cal{U})$
\item {\bf Axiom~6.5\ } $(\exists \cal{C})(\exists \cal{U})(\cof(\cal{C}) =
\cof(\fomom/\cal{U}))$ 
\end{itemize}\hyphenation{non-im-pli-ca-tion}
Various implications and non-implications between these axioms are
established in
\cite{nyik.c}. As well, it is observed that Axiom~2 is equivalent to
the equality $\frak{b} = \frak{d}$.
\section{Non-monotone  functions}
The definition of chains as $<^*$-increasing sequences of
nondecreasing functions in the 
axioms which appeared in \cite{nyik.c} may appear to be somewhat arbitrary
and one may wonder what results if
chains are defined differently. For the record, therefore,
the following definitions are offered.
\begin{defin}
If $\cal A$ is a subset of $\fomom$ then by
 a $(<^*,\cal{A})$-{\em chain}  will be meant a
subset of $\cal A$ which is well ordered by $<^*$.
By a $(\leq^*,\cal{A})${\em-chain}  will be meant a
subset of $\cal A$ which is well ordered by $\leq^*$.
For the purposes of this definition the most important subsets of
$\fomom$ are: the nonodecreasing functions, which will be denoted by
$\cal N$, the strictly increasing functions, which will be denoted by
$\cal S$, and $\fomom$.
\end{defin}
If   $x\in
\{<^*,\leq^*\}$ then Axiom~$N(x,\cal{A})$ will  denote the Axiom~$N$ 
with $\cal{C}$ being a
variable ranging over  $(x,\cal{A})$-chains ---
so Axiom~$N$ is the same as Axiom~$N(<^*,\cal{N})$. 

Fortunately, many of these axioms turn out to be equivalent and others
are simply false. The following simple observation of Nyikos
can be used to see this.
\begin{lemma}
    There is a mapping $\Psi : \fomom\rightarrow \fomom$ such that 
\label{integral} 
\begin{itemize}
\item $\Psi(f)$ is strictly increasing for every $f\in\fomom$
\item $f \leq \Psi(f)$ for every $f\in\fomom$
\item if $f\leq g$ and $f\neq^* g$ then $\Psi(f) <^* \Psi(g)$
\item if $f$ is nondecreasing and $f(n) < g(n)$ then $\Psi(f)(n) <
g(n)(n+1) + n$
\end{itemize}
\end{lemma}
\proof Define $\Psi(f)(n) = (\sum_{i=0}^nf(i)) + n$. It is easy to
check that $f \leq^*\Psi(f)$ and that $\Psi(f)$ is strictly increasing.
If $f\leq g$ and $f\neq^* g$ then there is 
some $m\in\omega$ such that $f(i) \leq
g(i)$ for all $i\geq m$. Since there are infinitely many $k\in
\omega$ such that $f(i)\neq g(i)$ is follows that there is some $M > m$
such that $\sum_{i=0}^Mg(i) > \sum_{i=0}^m(f(i)-g(i))$. Hence
$\Psi(f)(j) < \Psi(g)(j)$ for all $j \geq M$. 

Finally observe that if $f(n) < g(n)$ and $f$ is nondecreasing
then $f(i) < g(n)$ for each $ i \leq n$. Hence $\Psi(f)(n)\leq
f(n)(n+1) + n < g(n)(n+1) + n$. \stopproof

A consequence of Lemma \ref{integral} is that given any
$(\leq^*,\fomom)$-chain $\{f_{\xi} : \xi\in\lambda\}$ there is a
$(<^*,\cal{S})$-chain 
$\{f_{\xi}' : \xi\in\lambda\}$ such that $f_{\xi}\leq f'_{\xi}$ for each
$\xi\in\lambda$. Consequently,
 Axiom~$N(<^*,\cal{S})$,  Axiom~$N(<^*,\cal{N})$, 
Axiom~$N(<^*,\fomom)$, 
Axiom~$N(\leq^*,\fomom)$,  Axiom~$N(\leq^*,\cal{S})$ and Axiom~$N(\leq^*,
\cal{N})$ are all equivalent for $N\in \{2,4,6,6.5\}$. Therefore,
from now, if  $N\in \{2,4,6,6.5\}$, then Axiom~$N$ will be used to denote any and all of the Axioms
$N(x,\cal{A})$ where $x\in\{<^*,\leq^*\}$,
$\cal{A}\in\{\cal{N},\cal{S,\fomom}\}$.

Another consequence of Lemma \ref{integral} is that given any
$(\leq^*,\cal{N})$-chain $\{f_{\xi} : \xi\in\lambda\}$ there is a
$(<^*,\cal{S})$-chain 
$\{f_{\xi}' : \xi\in\lambda\}$ such that $f_{\xi}\leq f'_{\xi}$ for each
$\xi\in\lambda$ and, for each ultrafilter $\cal U$ and $g\in \fomom$,
the function $g$ is an upper bound for 
$\{f_{\xi} : \xi\in\lambda\}$ modulo $\cal U$ if and only if
$g\cdot (n + 1) + n$ is an upper bound for
$\{f_{\xi}' : \xi\in\lambda\}$ modulo $\cal U$. Consequently,
 Axiom~$N(<^*,\cal{S})$,  Axiom~$N(<^*,\cal{N})$, 
  Axiom~$N(\leq^*,\cal{S})$ amd Axiom~$N(\leq^*,
\cal{N})$ are all equivalent for $N\in \{1,3, 5\}$. 
Moreover, Axiom~$3(<^*,\fomom)$ and Axiom~$3(\leq^*,\fomom)$ are 
obviously both false because if $\cal U$
is any ultrafilter then it is possible to choose $X\subseteq \omega$
such that $X\notin \cal U$ and then find, using Lemma~\ref{integral},
 a $(<^*,\fomom)$-chain $\cal C$
such that $f(n) = 0$ for $f\in\cal C$ and $n\in X$. It follows from
this that Axioms~$1(<^*,\fomom)$ and $1(\leq^*,\fomom)$ are also false.
Therefore,
from now on Axiom~$N$ can be used to denote any and all of the Axioms
$N(x,\cal{A})$ where $x\in\{<^*,\leq^*\}$,
$\cal{A}\in\{\cal{N},\cal{S},\}$ and $N\in \{1,3\}$. Also
the notation Axiom~$5(\cal{N}) $ can be used to denote any and all of
 Axiom~$5(<^*,\cal{S})$,  Axiom~$5(<^*,\cal{N})$, 
  Axiom~$5(\leq^*,\cal{S})$ amd Axiom~$5(\leq^*,
\cal{N})$.

It is worth noting that Axiom~$5(\cal{N})$ is not equivalent to 
Axiom~${5(<^*,\fomom)}$ or Axiom~${5(\leq^*,\fomom)}$.
The reason for this is that it will been shown,  in 
Theorem \ref{2implies5}, that
Axiom~2 implies Axiom~$5(\cal{N})$; but the same is not
true of
Axiom~$5(<^*,\fomom)$ since the next result shows that Axiom~$5(<^*,\fomom)$
 fails assuming
$2^{\aleph_0 } = \aleph_1$. It is obvious that Axiom~2 holds if
$2^{\aleph_0 } = \aleph_1$. The following definition
 will be used to establish this and appears to be
central in the context of non-monotone functions.
\begin{defin}
 If $\cal{A}\subseteq\fomom$ then define $\cal{I}_b(\cal{A})$ to be the
set of all $X\subseteq\omega$ such that $\{f\restriction X : f\in
\cal{A}\}$ is bounded.
\end{defin}
Notice that $\cal{I}_b(\cal{A})$ is an ideal and that
$\cal{I}_b(\cal{A})$ is proper if and only if $\cal{A}$ is an 
unbounded subset of $\fomom$.
It is also worth observing that if $\cal{C}\subseteq\fomom$ and $\cal{U}$
is an ultrafilter and
$\cal{C}$ is unbounded modulo $\cal{U}$ then
$\cal{U}\cap\cal{I}_b(\cal{C}) = \emptyset$. 
\begin{lemma}
    If  \label{CHlem}  there is a sequence $\{X_{\xi} : \xi\in\omega_1\}$
 of subsets of $\omega$ such that
\begin{itemize}
    \item $X_\xi\subseteq^* X_\eta$ if $\xi\in\eta$
\item  $X_\eta\setminus X_\xi$ is infinite if $\xi\in\eta$
\item there exists a family $\{g_\xi : \xi\in\omega_1\}\subseteq\fomom$
such that for every $f\in\fomom$ there is $\xi\in\omega_1$ such that
$g_\xi\restriction X_{\xi + 1}\setminus X_\xi \not\leq^* f\restriction
 X_{\xi + 1}\setminus X_\xi$
\end{itemize} then
there is an unbounded $(<^*,\fomom)$-chain,  $\cal C$, such that
$\cal{I}_b(\cal{C})$ contains
$ \{\{n\in\omega : f(n) \geq n\} : f\in\cal{C}\}$.
\end{lemma}
\proof
Let $\{X_{\xi} : \xi\in\omega_1\}$ and $\{g_{\xi} :
\xi\in\omega_1\}$ satisfy the hypothesis of the lemma and, without
loss of generality, assume that $g_\xi(n)\geq n$ for all $\xi$ and $n$. 
Let $\{h_{\xi} : \xi\in
\omega_1\}$ be a $<^*$-increasing sequence of functions such that
$h_{\xi}(n) < n$ for all $n\in\omega$. A standard induction argument can
now be used to construct $\{f_{\xi} : \xi\in\omega_1\}$ such that
\begin{itemize}
\item $f_{\xi}\restriction \omega\setminus X_{\xi + 1} = 
h_{\xi}\restriction \omega\setminus X_{\xi + 1}$
\item if $\xi\in\eta$ then $f_{\xi}\leq^* f_{\eta}$
\item if $\xi\in\eta$ then $f_{\xi}\restriction  X_{\xi}
=^* f_{\eta}\restriction  X_{\xi}$
\item $f_{\xi}\restriction  X_{\xi + 1}\setminus X_\xi = g_{\xi}
\restriction  X_{\xi + 1}\setminus X_\xi$
\end{itemize} and this clearly suffices.
\stopproof

Notice that $\omega_1$ is crucial to the proof of Lemma~\ref{CHlem}
and can not be replaced by a larger cardinal. The reason is that
the inductive construction relies on the fact that if $\{f_n : n\in\omega\}
$ is a family of partial functions from $\omega$ to $\omega$  such that 
$f_n =^* f_{n+1}\restriction \dom(f_n)$ then, there is a single function
$f$ such that $f_n\subseteq^* f$ for each $n\in \omega$. A Hausdroff gap type of construction shows that this is not possible if
 $\omega_1$ is replaced by some larger cardinal.
It is for the same reason that $\omega_1$ appears in the next corollary.
\begin{corol}
    If $\frak{d}=\aleph_1$ then \label{CHlemcor}
there is an unbounded $(<^*,\fomom)$-chain,  $\cal C$, such that
$\cal{I}_b(\cal{C})$ contains
$ \{\{n\in\omega : f(n) \geq n\} : f\in\cal{C}\}$.
\end{corol}
\proof Let 
 $\{X_{\xi} : \xi\in\omega_1\}$ be any sequence
 of subsets of $\omega$ such that
 $X_\xi\subseteq^* X_\eta$ 
  $X_\eta\setminus X_\xi$ is infinite if $\xi\in\eta$.
Then $\{g_\xi : \xi\in\omega_1\}$ can be any dominating family.
\stopproof

The next result shows that
	 if $\frak{d} = \aleph_1$ then Axiom~$5(<^*,\fomom)$ fails.
\begin{theor}
    If $\frak{d} = \aleph_1$ then there is an unbounded
$(<^*,\fomom)$-chain 
 which is \label{CHth} bounded modulo any ultrafilter.
\end{theor}
\proof  Use 
Corollary~\ref{CHlemcor} to find
an unbounded chain $\cal{C}\subseteq\fomom$ such that
$\cal{I}_b(\cal{C})$ contains 
$ \{\{n\in\omega : f(n) \geq n\} : f\in\cal{C}\}$.
If $\cal{V}$ is any ultrafilter such that $\cal C$ is unbounded modulo
$\cal V$ then it must be that
$\cal{V}\supseteq
 \{\{n\in\omega : f(n) \geq n\} : f\in\cal{C}\}=\emptyset$. Then,
it is clear that the identity function is an upper bound for $\cal{C}$.
\stopproof

It seems that Axiom~$5(<^*,\fomom)$ is very strong and
Axiom~${5(\leq^*,\fomom)}$ is potentially even stronger.
Nevertheless, Axiom~${5(\leq^*,\fomom)}$ is consistent and does not imply
Axiom~1. This is implied by the next sequence of results. The 
question of which  of the axioms are implied by Axiom~$5(\leq^*,\fomom)$ 
is mostly open however.

    The Open Colouring Axiom was first considered
by Abraham, Rubin and Shelah in \cite{ab.ru.sh} and later
strengthened by Todorcevic \cite{todo.pp}.
\begin{defin}
     The Open Colouring Axiom states that if $X\subseteq \RR$ and
 $\cal{V}\subset [X]^2$
is an open set\footnote{Here $[X]^2$ can be thought of as the set of points in
$X^2$ above the diagonal.}
 then either there is $Y\in [X]^{\aleph_1}$ such that $[Y]^2 \subset
\cal{V}$ or there exists a partition of $X=\cup_{n\in\omega}X_n$
such that
$[X_n]^2\cap\cal{V}=\emptyset$ for each $n\in\omega$.
$\RR$ can be replaced by any second countable
space in the statement of the Open Colouring Axiom.
\end{defin}

\begin{theor}
If the Open Colouring Axiom holds and\label{oca}
 $\{(h_{\alpha},g_{\alpha}) : \alpha\in\lambda\}$ satisfies
\begin{itemize}
\item $\lambda$ is a regular cardinal greater than $\omega_1$ 
\item $\dom(h_\xi) = \dom(g_\xi) = X_\xi$ for $\xi\in\lambda$
\item if $\xi\in \eta$ then $X_\xi\subseteq^* X_\eta$
\item if $n\in X_\xi$ then $g_\xi(n)\leq h_\xi(n)$
\item if $\xi \in\eta$ then $g_\xi\leq^*g_\eta\restriction X_\xi \leq^* h_\xi$
\end{itemize}
then  there exists a function $f:\omega\rightarrow\omega$ such that
$g_\xi \leq^* f\restriction X_\xi$ for all $\xi\in\lambda$.
\end{theor}
\proof
To begin, identify $\lambda$ with the subspace of  the reals
$\{(h_\xi,g_\xi) : \xi\in \lambda\}$ --- the reals are being considered
as $(\fomom)^2$ or, in other words, the irrationals. Define $$V = 
\{\{\alpha,\beta\} \in [\lambda]^2 : \alpha\in\beta\AND (\exists n)(g_\beta(n)
> h_\alpha(n))\}$$
and observe that $V$ is open. From the Open Colouring Axiom it follows that
there are only two possibilities.

The first is that there is a partition $\lambda = \cup_{n\in\omega}X_n$
such that $[X]^2\cap V=\emptyset$ for each $n\in\omega$. In this case there
must be some $n\in\omega$ such that $X_n$ is cofinal in $\lambda$. Choose
$\xi\in\lambda$ such that $\cup_{\eta\in\xi}X_\eta  = \cup_{\eta\in\lambda}X_\eta $ and let $f(n) = \min\{h_\eta(n) : \eta\in\xi\}$. Now, if $\beta > \xi$
and $n\in X_\beta$ then there is some $\eta\in\xi$ such that $n\in X_\eta$
and hence $n\in\dom(f)$. Moreover, if $\eta\in\xi$ and $n\in X_\eta$ then
 $g_\beta(n)\leq h_\eta(n)$ and so $g_\beta(n) \leq f(n)$.
So $f$ is the desired function.

The second possibility is that there is $X\in [\lambda]^{\aleph_1}$ such that
$[X]^2\subset V$. Since $\lambda \geq \omega_2$ it is possible to choose
some $\xi\in\lambda$ such that $X\subseteq \xi$. It is then possible to choose
$M\in \omega$, $g:M\rightarrow \omega$ and $Y\in [X]^{\aleph_1}$ such that
\begin{itemize}
\item if $\eta\in Y$ then $X_\eta \setminus M\subseteq X_\xi$
\item if $\eta\in Y$ and $n\in X_\eta \setminus M$ then $g_\eta(n)\leq g_\xi(n)
\leq h_\eta(n)$
\item $g_\eta\restriction M = g$
\end{itemize}
Then if $\mu\in\eta$ and $\{\mu,\eta\}\in [Y]^2$ and $n\in X_\mu\cap X_\eta$
it must be that either $n \geq M$ or $n < M$. In the first case it follows that
$n\in X_\xi$ and so $g_\eta(n)\leq g_\xi(n)\leq h_\mu(n)$. In the second
case it may be concluded  that $g_\eta(n) = g(n) = g_\mu(n) \leq h_\mu(n)$. It follows
that $\{\mu,\eta\}\notin V$ which is a contradiction.
\stopproof

\begin{theor}
The conjunction of Axiom~2\label{pfaimp5} and the Open Colouring
Axiom implies Axiom~5$(\leq^*,\fomom)$.
\end{theor}
\proof To begin, recall that it was shown in \cite{nyik.c} that
Axiom~2 implies that
$\frak{b} = \frak{d}$. Hence it is possible to choose
a  $(\leq^*,\fomom)$-chain $ \{d_\xi : \xi \in \frak{d}\}$
which is also a dominating family in $\fomom$. Also,
if $\cal{C}$ is any  $(\leq^*,\fomom)$-chain then $\cal C$
is of the form $ \{g_\xi : \xi \in \frak{d}\}$.
Define $E(\eta,\xi) = \{n\in \omega : g_\xi(n) \geq d_\eta(n)\}$.

Next, let $\{\frak{M}_\xi : \xi\in \frak{d}\}$ be a sequence of elementary
submodels of $(H(\frak{c}^+),\in)$ such that
\begin{itemize}
\item $\card{\frak{M}_\xi} < \frak{d}$ for each $\xi \in \frak{d}$
\item $ \{g_\xi : \xi \in \frak{d}\} \in \frak{M}_\eta$
and $ \{d_\xi : \xi \in \frak{d}\}\in \frak{M}_\eta$ 
for each $\eta \in \frak{d}$
\item $\eta \in \frak{M}_\eta$
\item  $\frak{M}_\eta\in \frak{M}_\zeta$ 
for each $\eta \in \zeta \in \frak{d}$
\end{itemize}
and let $\mu(\xi) = \frak{M}_\xi\cap\frak{d}$.
Define $\cal F$ to 
 be the filter generated by $$\{E(\mu(\xi),\mu(\xi + 1)) :
\xi \in\frak{d}\AND \xi\mbox{ is odd }\}$$
and observe that if $\cal F$  is a proper filter then $\cal C$ will be 
cofinal in $\fomom$ modulo $\cal U$ for any ultrafilter extending $\cal F$.

Hence it suffices to show that $\cal F$ is proper. To this end let
 ${\cal F}_\rho$
 be the filter generated by $$\{E(\mu(\xi),\mu(\xi + 1)) :
\xi \in \rho\AND \xi\mbox{ is odd }\}$$
and prove by induction that each  ${\cal F}_\rho$ is proper.
Moreover, it will be shown by induction that  ${\cal F}_\rho\cap
\cal{I}_b(\cal{C}) = \emptyset$.
If $\rho  = 0 $, $\rho$ is odd
 or $\rho$ is a limit  then there is nothing to do so
suppose that $\rho = \rho'+ 1$, where $\rho'$ is odd, and that  
${\cal F}_{\rho'}$ is a proper filter
such that   ${\cal F}_{\rho'}\cap
\cal{I}_b(\cal{C}) = \emptyset$. 

Notice that ${\cal F}_{\rho'}\in \frak{M}_{\rho'}$ because $\rho'$ 
is odd. 
Therefore it suffices to show that for each $B\in\cal{I}_b(\cal{C})^+$ there
is some $\theta\in \frak{d}$ such that $E(\mu(\rho'),\theta)\cap B\in
\cal{I}_b(\cal{C})^+$ --- the reason being that the elementarity of
$\frak{M}_{\rho' + 1}$ will guarantee that
$E(\mu(\rho'),\mu(\rho'+1))\cap B\in
\cal{I}_b(\cal{C})^+$ for each $ B\in
\cal{I}_b(\cal{C})^+$. Elementarity also assures that it may as well be 
assumed that $B\in \frak{M}_{\rho'}$.
But if there is some $ B\in
\cal{I}_b(\cal{C})^+$ such that $E(\mu(\rho'),\theta)\cap B\in
\cal{I}_b(\cal{C})$ for each $\theta \in \frak{c}$ then it is possible to find
$h_\theta$ such that 
\begin{itemize}
\item $\dom(h_\theta) = B\cap 
E(\mu(\rho'),\theta)$ for each $\theta\in\frak{d}$
\item $h_\theta(n)\geq g_\theta(n)$ for every $n\in
B\cap E(\mu(\rho'),\theta)$ and for each $\theta\in\frak{d}$
\item $g_\xi\restriction \dom(h_\theta) \leq^* h_\theta$ 
for each $\theta\in
\xi\in\frak{c}$
\end{itemize}
It follows that $\{(h_\theta,g_\theta\restriction B\cap E(\mu(\rho'),\theta)) :
\theta \in \frak{d}\}$ satisfies the hypothesis of Lemma~\ref{oca}.
Since the Open Colouring Axiom is being assumed,
 there is a function $f\in\fomom$ such that
$g_\theta\leq^* f\restriction B\cap E(\mu(\rho'),\theta)$ for each $\theta \in
\frak{d}$. It follows that for each $\theta \in \frak{d}$ there 
are only finitely many $n\in B$ such that $g_\theta(n) > \max\{d_{\rho'}(n),f(n)\}$
contradicting that $B\notin \cal{I}_b(\cal{C})$.
\stopproof

Notice that it is shown in \cite{todo.pp} that the Proper Forcing Axiom
implies the hypothesis of Theorem~\ref{pfaimp5}. Moreover it is a
Corollary that Axiom~${5(\leq^*,\fomom)}$ does not imply Axiom~1 because
it is easy to check that Martins' Axiom --- and hence the Proper Forcing
Axiom --- implies that Axiom~1 fails. In particular, it is possible to inductively define
a $(<^*,\cal{S})$-chain  no member of which dominates the
exponential\footnote{The exponential function is not crucial here but
some quickly growing function must be used. For example, although the
identity function is strictly increasing it can not be used because it
is the minimal strictly increasing function.}  
function.

It has already been mentioned that the next lemma can be used to show
that Axiom~$5(\cal{N})$ is not equivalent 
to Axiom~$5(\leq^*,\fomom)$ or Axiom~${5(<^*,\fomom)}$. It
 will also be used in the proof of Theorem \ref{6not5.5} but also
has some interest on its own since it provides a sufficient
condition for Axiom~$5(\cal{N})$ to hold. Thus, it 
will be used to show that
Axiom~2 implies Axiom~$5(\cal{N})$.
\begin{lemma}
    If  $\frak{d}$  is regular 
and\label{suff6}
 $\{c_{\xi} : \xi\in
\frak{d}\}\subseteq \fomom$ is $\leq^*$ increasing and, moreover,
$\{c_{\xi}\restriction A : \xi\in
\frak{d}\}$ is unbounded in $\presup{A}{\omega}$ for each $A\in
[\omega]^{\aleph_0}$ then there is an ultrafilter $\cal{U}$ such that
$\{c_{\xi} : \xi\in
\frak{d}\}$ is cofinal in  $\fomom/\cal{U}$.
\end{lemma}
\proof
Let $\cal{D}\subseteq \fomom$ be a cofinal family in $\fomom$ of cardinality
$\frak{d}$.
Let $\{\frak{M}_{\xi} : \xi\in\frak{d}\}$ be an increasing sequence
of elementary submodels of $$(H(\frak{c}^+), \{c_{\xi} : \xi\in
\frak{d}\},\cal{D},\in)$$ such that $\frak{M}_{\xi}\cap\frak{d} =
\alpha(\xi)\in\frak{c}$
for each $\xi\in\frak{c}$
 and $\cup_{\xi\in\frak{d}}\frak{M}_{\xi}\supseteq
\cal{D}$
 --- this is possible because $\frak{d}$ is regular. Let $$\cal{U} =
\{\{n\in\omega : f(n) \leq c_{\alpha(\xi)}(n)\} : \xi\in\frak{d}\AND
f\in\frak{M}_{\xi}\}$$ and note that it suffices to show that
 this is a base for a filter.

That $\cal{U}$ has the finite intersection property can be
established by induction. Let
$B(\xi,f) = \{n\in\omega : f(n) \leq c_{\alpha(\xi)}(n)\}$
for $\xi\in\frak{c}$ and 
$f\in\frak{M}_{\xi}$ and suppose that
$\card{\cap\cal{A}} = \aleph_0$ for each $\cal{A}\in [\cal{U}]^m$
--- the case $m = 1$ is an easy consequence of elementarity.
Now let
$$\{B(\xi_0,f_0),B(\xi_1,f_1),\dots,B(\xi_m,f_{m})\}\in [\cal{U}]^{m+1}$$
 and suppose that $\xi_i \leq \xi_{i+1}$ for each $i$. If $\xi_{m-1}
= \xi_m$ then $\{f_{m-1},f_m\}\subseteq\frak{M}_{\xi_m}$ and so the
elementarity of $\frak{M}_{\xi_m}$ ensures that there is some
$g\in\frak{M}_{\xi_m}$ such that $f_{m-j}\leq^* g$ for each $j\in
2$. Hence
$B(\xi_0,f_0)\cap B(\xi_1,f_1)\dots B(\xi_m,f_{m})$ contains
$$B(\xi_0,f_0)\cap B(\xi_1,f_1)\cap \dots\cap B(\xi_{m-2},f_{m-2})\cap
B(\xi_m,g)$$ and this set is
 infinite by the induction hypothesis.

On the other hand,
if $\xi_{m-1} \in \xi_m$ then
$$B = B(\xi_0,f_0)\cap B(\xi_1,f_1)\cap\dots\cap B(\xi_{m-1},f_{m-1})$$
is infinite by the induction hypothesis and, moreover, $B$ belongs to
$\frak{M}_{\xi_{m}}$ because all the parameters defining it do.
Since
$\{c_{\xi}\restriction B : \xi\in
\frak{d}\}$ is unbounded in $\presup{B}{\omega}$ it follows that
there must be some $\mu\in \frak{M}_{\xi_m}$ such that
$f_{m}\restriction B \not\leq^*
c_{\mu}\restriction B$ and so  $f_m\restriction B \not\leq^*
c_{\alpha(\xi)}\restriction B$. Since
$B(\xi_0,f_0)\cap B(\xi_1,f_1)\cap\dots\cap B(\xi_{m},f_{m}) =
B\cap B(\xi_m,f_m)$ this is enough.
\stopproof
\begin{theor}
    Axiom~2 implies \label{2implies5} Axiom~$5(\cal{N})$.
\end{theor}
\proof In \cite{nyik.c} it is shown that Axiom~2 is equivalent to
the equality $\frak{b} = \frak{d}$. Since $\frak{b}$ is regular it
follows that $\frak{d}$ is regular. Moreover, if
$\cal{C}\subseteq\fomom$ is an unbounded $(\leq^*,\cal{N})$-chain 
then $\cof{(\cal{C})} =\frak{d}$.
Since $\cal{C}$ consists of nondecreasing functions it is clear that
$\{c\restriction A : c\in\cal{C}\}$ is unbounded for each infinite
set $A$. Hence, by Lemma \ref{suff6}, it follows that there is an ultrafilter
$\cal{U}$ such that $\cal{C}$ is unbounded modulo $\cal{U}$.
\stopproof

\section{Oracle Chain Conditions and Locally Cohen Partial Orders}
It will be shown that there is a model of set theory where Axiom~6.5
fails. This answers the first two questions in Problem 5 of \cite{nyik.c}.
C. Laflamme has remarked that in some models of NCF (see  \cite{blass.appl}
for an overview of this area) Axiom 6.5 fails as well because it is possible to provide a classification of chains  in these models.
The restriction to chains does not play an
important role in this theorem and, in fact, the theorem is 
slightly stronger than required --- at least formally --- because of this.
\begin{theor}
    There is a model \label{not6} where $\cof(\fomom/\cal{U}) =
\omega_2$ for every ultrafilter $\cal{U}$ 
but every unbounded subset of $\fomom$ has an unbounded subset of
size  $\aleph_1$.
\end{theor}
\proof 
The plan of the proof is to start with a model $V$ in which
$\lozenge^*_{\omega_1}$ and $\lozenge_{\omega_2}(\omega_1)$
--- in other words, the trapping of
subsets of $\omega_2$ occurs
at ordinals of cofinality $\omega_1$ in $\omega_2$ --- both hold. 
In this model a finite support iteration
$\{(\PP_{\xi}, \QQ_{\xi}) : \xi\in\omega_2\}$ will be constructed
along with a sequence of oracles \cite{shel.pf}
 $\{\frak{M}_{\xi} : \xi\in\omega_2\}$ --- more precisely,
$\frak{M}_{\xi }$ is a $\PP_{\xi}$-name for an oracle. 
The oracles will be chosen
so that if $\{g_{\eta} : \eta\in\omega_1\}$ is a $\PP_{\xi}$-name,
guessed by the  $\lozenge_{\omega_2}(\omega_1)$ sequence,
for an unbounded subset of $\fomom$ then
$\frak{M}_{\xi }$ is chosen so that if $\QQ$ is any partial order
satisfying the $\frak{M}_{\xi }$-chain condition then forcing
with $\QQ$ does not destroy the unboundedness of  $\{g_{\eta} :
\eta\in\omega_1\}$. Provided that $\PP_{\omega_2}/\PP_{\xi}$
satisfies the $\frak{M}_{\xi}$-chain condition, it will follow
that every  unbounded subset of $\fomom$ has cofinality $\omega_1$
because every unbounded subset is reflected at some initial stage by the
 $\lozenge_{\omega_2}(\omega_1)$ sequence. The rest of the result will
follow once it is shown how to construct $\QQ_{\xi}$ satisfying the
$\frak{M}_{\xi}$-chain condition and adding an upper bound to
any given sequence from some ultrapower of the integers.

The construction of $\{(\PP_{\xi}, \QQ_{\xi}) : \xi\in\omega_2\}$
and $\{\frak{M}_{\xi} : \xi\in\omega_2\}$ is, of course, done by induction.
If $\xi$ is a limit then $\PP_{\xi}$ is simply the direct limit of
$\{\PP_{\mu} : \mu\in\xi\}$. The construction of $\frak{M}_{\xi}$
and $\QQ_{\xi}$
does not depend on whether or not $\xi$ is a limit.

Given $\PP_{\xi}$, use
the results of pages  124 to 127 of \cite{shel.pf} to find a
$\PP_{\xi}$-name for a single oracle $\frak{N}_{\xi}$ such that if
$\QQ$ satifies the $\frak{N}$-chain condition then it satisfies the
$\frak{M}_{\mu}$ chain condition for each $\mu\in \xi $. Let
$C_{\xi}$ be the set guessed by the $\lozenge_{\omega_2}(\omega_1)$ 
sequence at $\xi$. If $C_{\xi}$ is not a $\PP_{\xi}$-name 
for an unbounded subset of $\fomom$ then let $\frak{M}_{\xi } =
\frak{N}_{\xi}$. Otherwise, use Lemma 2.1 on page 122 of
\cite{shel.pf} to find an oracle $\frak{M}$ such that if $\QQ$
satisfies the $\frak{M}$-chain condition then the subset
 $C_{\xi}\subseteq\fomom$
remains unbounded after forcing with $\QQ$. The use of Lemma 2.1
requires checking that if $C\subseteq\fomom$ is an unbounded chain then
adding a Cohen real will not destroy its unboundedness. This is a
result of the folklore which can be found in \cite{step.20}.
Then use the results of pages  124 to 127 of \cite{shel.pf} to find a
single oracle $\frak{M}_{\xi}$ such that any $\QQ$ which satisfies
the $\frak{M}_{\xi}$-chain condition will also satisfy the $\frak{M}$-chain
condition and the $\frak{N}_{\xi}$-chain condition.

Suppose that the $\lozenge_{\omega_2}(\omega_1)$ sequence has also trapped a
filter $\cal{U}_{\xi}$ --- which is an ultrafilter in the
intermediate generic extension  by $\PP_{\xi}$ --- 
and an increasing sequence $\{f_{\mu}^{\xi} :
\mu \in \omega_1\}$ in the  reduced 
power of the integers modulo $\cal{U}_{\xi}$.
(So it is being assumed that, by some coding, the
$\lozenge_{\omega_2}(\omega_1)$ sequence traps triples of sets --- the first
component of the triple at $\xi$ is a candidate for $C_{\xi}$ in the
construction of $\frak{M}_{\xi}$ while the second and third
components are candidates for the ultrafilter $\cal{U}_{\xi}$ and
the sequence $\{f_{\mu}^\xi : \mu\in\omega_1\}$.) 
The only thing left to do is to construct $\QQ_{\xi}$ satisfying the
$\frak{M}_{\xi}$-chain condition and adding an upper bound for 
$\{f_{\mu}^\xi : \mu\in\omega_1\}$ in the reduced power modulo
$\cal{U}_{\xi}$. 

Let $\frak{M}_{\xi} = \{M_{\xi}^{\mu} : \mu\in\omega_1\}$. 
The partial order $\QQ_{\xi}$ is constructed by induction on
$\omega_1$ in $V^{\PP_{\xi}}$ --- it will be similar to the forcing
which adds a dominating real but with extra side conditions. 
In particular, a sequence of partial functions $\{S_{\mu}
: \mu\in\omega_1\}\subseteq\fomom$ is constructed by induction 
on $\omega_1$ and
$\QQ_{\xi}^{\mu}$ is defined to be the set of all pairs
$(F,\Gamma)$ such that $F:k\rightarrow\omega$ is a finite partial function
and  $\Gamma\in [\mu]^{<\alpeh_0}$.
The ordering on $\QQ_{\xi}^{\mu}$ is defined by $(F,\Gamma) \leq
(F',\Gamma')$ provided that $\Gamma\subseteq \Gamma'$,
$F\subseteq F'$ and $F'(j)\geq S_{\gamma}(j)$ for
$\gamma\in\Gamma$ and  $j\in(\dom(S_{\gamma})\setminus\dom(F))$. 
Moreover, the functions $S_{\mu}$ will be constructed so that
$\dom(S_{\mu})\in \cal{U}_{\xi}$ and $S_{\mu}(j)\geq f_{\mu}^{\xi}(j)$
for each $j\in \dom(S_{\mu})$.
It is easy to see that $\{(F,\Gamma) : \mu\in\Gamma\}$ is dense in
$\QQ_{\xi}^{\eta}$ for every $\mu\in\eta$ and so if $G$ is 
$\QQ_{\xi} = \QQ_{\xi}^{\omega_1}$ generic over $V^{\PP_{\xi}}$ then
$\cup\{F : (\exists \Gamma)((F,\Gamma)\in G)\}$ is an upper bound
for $\{f_{\mu}^\xi : \mu\in\omega_1\}$ in the  reduced power with respect to
$\cal{U}_{\xi}$. It therefore
suffices to construct $\{S_{\mu} : \mu\in\omega_1\}$ so that for
every $\mu\in\omega_1$,  every dense open subset of  $\QQ_{\xi}^{\mu}$
which belongs to $M_{\xi}^{\mu}$ remains predense in $\QQ_{\xi}^{\mu+1}$.
This, of course, will ensure that $\QQ_{\xi} = \QQ_{\xi}^{\omega_1}$ satisfies
the $\frak{M}_{\xi}$-chain condition.

Suppose that $\{S_{\mu} : \mu\in\eta\}$ have been constructed. Let
$\cal{A}$ be the set of all of the dense open subsets of
$\QQ_{\xi}^{\eta}$ which belong to  $M_{\eta}$ --- this includes
all those 
 dense open subsets of
$\QQ_{\xi}^{\eta}$ which belong to
$M_{\zeta}$ some $\zeta\in
\eta $. Choose $h\in\fomom$ to be some function
which dominates all members of $M_{\eta}$; in other words, if
$g\in\fomom\cap M_{\eta}$ then $g\leq^* h$. Let 
$\{(A_i,(F_i,\Gamma_i)) : i\in\omega\}$ enumerate
$\cal{A}\times\QQ_{\xi}^{\eta}$. Now choose, by induction on
$\omega$, integers $\{K_i : i\in\omega\}$ such that
$K_i < K_{i+1}$ and $K_0 = 0$. Given $K_i$, define
$\bar{F_j^i}\supset F_j$ for each $j\leq
i$ such that if $\dom(F_j)\subseteq K_i$ then $\dom(\bar{F_j^i}) = K_i$,
 $(F_j,\Gamma_j)\leq  (\bar{F_j^i},\Gamma_j)$ and $\bar{F_j^i}(k)\geq
h(k)$ if $k\in\dom({\bar{F_j^i}\setminus F_j})$. Now choose
$(F^i_j,\Gamma_j^i)\in A_j$ such that $(F^i_j,\Gamma_j^i)\geq
(\bar{F_j^i},\Gamma)$. Let $K_{i+1}$ be such that $\dom(F_j^i)\subset
K_{i+1}$ for each $j \leq i$. 

Let $X_m = \cup_{i\in\omega}[K_{2i+m},K_{2i+m+1})$ for $m\in 2$ and
note that there exists $m'\in 2$ such that $X_{m'}\in \cal{U}_{\xi}$.
Let $S_{\mu} = h\restricts X_{m'}$; the reason being that $\cal{U}_{\xi}
$ is an ultrafilter in $V^{\PP_{\xi}}$ and $X_0\in V^{\PP_{\xi}}$.
 To see that this 
definition of $S_{\mu}$ ensures
that every dense open subset of $\QQ_{\xi}^{\mu}$
which belongs to $M_{\xi}^{\mu}$ remains predense in
$\QQ_{\xi}^{\mu+1}$ let $(F,\Gamma) \in \QQ_{\xi}^{\mu+1}$ and let
$D\in M_{\mu}$ be dense open in $\QQ_{\xi}^{\mu}$. It
follows that $(F,\Gamma\setminus\{\mu\}) \in \QQ_{\xi}^{\mu}$. To simplify
notation assume that $m = 1$. 
Choose $j$ such that $\dom (F)\subseteq K_{2j}$
and such that
$(D,(F,\Gamma\setminus \{\mu\}) = (A_k,(F_k,\Gamma_k))$
 for some $k\leq 2j$. \hyphenation{since}
Since $\dom({F_k^{2j}})\subseteq K_{2j}$
it follows that $\bar{F_k^{2j}}(n)\geq h(n)$
if $n\in\dom({\bar{F_k^{2j}}\setminus F_k})$. Moreover $[K_{2j},
K_{2j+1})\cap \dom S_{\mu} = \emptyset$. Hence
$(F_k^{2j},\Gamma_k^{2j})$ is compatible with $(F,\Gamma)$.
\stopproof

The methods of the previous theorem can also be used to show that it
is consistent that Axiom~6 holds but Axiom~5.5 fails. 
In establishing this it will be helpful to  introduce the following definition.
\begin {defin}
    A partial order $(\PP, \leq)$ will called {\em locally Cohen} if for every
$X\in [\PP]^{\aleph_0}$ there is $Y\in [\PP]^{\aleph_0}$ such that
$X\subseteq Y$ and $Y$ is completely embedded in $\PP$ --- in other
words, if $A\subseteq Y$ is a maximal antichain in the partial order
$(Y,\leq\cap Y\times Y)$ then it is also maximal in $(\PP,\leq)$.
\end{defin}
The notion of locally Cohen partial orders has already been isolated and
investigated by W. Just in \cite{just} who refers to locally Cohen partial
orders as {\em harmless}. 
The motivation of Just was that any locally Cohen forcing satisfies the 
oracle chain condition  for every oracle.

Let $\SS(\lambda)$ be the canonical partial order for adding a scale
of length $\lambda$ in $\fomom$ with finite conditions. To be precise,
a condition $p$ belongs to $\SS(\lambda)$ if and only if $p:\Gamma_p\times n_p
\rightarrow \omega$ is a
 function and $\Gamma_p\in [\lambda]^{<\aleph_0}$ and $n_p\in\omega$. 
The ordering on $\SS(\lambda)$ is $\leq$ defined by $p \leq q$ 
if and only if:
\begin{itemize}
    \item $p \subseteq q$
\item if $\{\xi , \eta\} \subseteq\Gamma_p$ and
$\xi\in \eta$ then $q(\eta,m) \geq q(\xi,m)$ 
for every $m\in n_q\setminus n_p$
    \end{itemize}

It should be noted that that $\SS(\lambda)$ is also the finite 
support iteration
of length $\lambda $ of the partial orders $\{\DD(\xi) : \xi\in \lambda\}$ 
where
$\DD(\xi)$ is the finite condition forcing for adding a
nondecreasing function  ---
which will be denoted by $c_{\xi}$ --- which dominates all the reals
$\{c_{\eta} : \eta \in \xi\}$.
\begin{lemma}
    For any ordinal $\lambda$ the partial order $\SS(\lambda)$
is locally Cohen.
\end{lemma}
\proof Given
$X\in [\PP]^{\aleph_0}$ let $Y\in [\lambda]^{\aleph_0}$ be any 
set such that $\SS(Y)\supseteq X$ --- $\SS(Y)$ can be defined for any set of
ordinals in the same way that $\SS(\lambda)$ is defined 
for an ordinal $\lambda$. To see that $\SS(Y)$ is completely embedded in
$\SS(\lambda)$ let $A \subseteq \SS(Y)$ be a maximal antichain in $\SS(Y)$.
If $p \in \SS(\lambda)$ then let $\Gamma' = 
\Gamma\cap Y$ and $p' = p\restriction \Gamma'\times n_p$. Since $p'\in
\SS(Y)$ there must be some $q\in A$ such that $q' =p'\cup q\in \SS(Y)$ and
$p \leq q'$. Define $q'':(\Gamma_q\cup \Gamma_p)\times n_q \rightarrow \omega$
by $$q''(\xi,j) = \left\{\begin{array}{ll}
			    q'(\xi,j)  & \IF \xi\in Y  \\
			    \max\{q'(\eta,j) : \eta \in Y\cap\xi\}
& \IF j\notin n_p \AND \xi\notin Y\\
p(\xi,j) & \IF\xi\notin Y \AND j\in n_p
			 \end{array}\right.$$
It is easy to check that $p \leq q''$.
\stopproof

\begin{theor}
\label{6not5.5}    There is a model of set theory where 
\begin{itemize}
    \item $2^{\aleph_0} = \aleph_2$\item $\frak{b} =
\aleph_1$
\item there is a an unbounded $(\leq^*,\fomom)$-chain of length $\omega_2$
\item the cofinality of any ultrapower of the integers is $\omega_2$
\end{itemize}
\end{theor}
\proof
As in the proof of Theorem \ref{not6}, let $V$ be a model 
of set theory where $\lozenge_{\omega_1}^*$ and
$\lozenge_{\omega_2}(\omega_1)$ 
are both satisfied.
Let $H$ be $\SS(\omega_2)$ generic over $V$ and define $H_{\xi} =
H\cap\SS(\xi)$. Let $c_{\xi}(n) = (\cup H)(\xi,n)$ and observe that
$\{c_{\xi} : \xi\in\omega_2\}\subseteq\fomom$ is increasing with
respect to $\leq^*$ and
$\{c_{\xi} : \xi\in\omega_2\}$ is not bounded. 
If  an oracle chain condition forcing extension of $V[H]$ can be found
which preserves the unboundedness of
$\{c_{\xi} : \xi\in\omega_2\}$ and in which
the cofinality of any ultrapower of the integers is $\omega_2$
then the result will follow because $\frak{b} = \aleph_1$ is easily
preserved by the oracle chain condition.

To do this,
 construct  $\{(\PP_{\xi}, \QQ_{\xi}) : \xi\in\omega_2\}$
and
$\{\frak{M}_{\xi} : \xi\in\omega_2\}$ exactly as in the proof of Theorem
\ref{not6}
except that
$\frak{M}_{\xi}$ is chosen to be an oracle in the model $V[G,H_{\xi}]$ where
$G$ is generic
over
$\PP_{\xi}$.
There is no
problem in
doing this
because
$\SS(\xi)$ is
locally Cohen
and hence
satisfies the
$\frak{M}_{\mu}$ chain condition for each $\mu\in\xi$ --- indeed, $\SS(\xi)$ 
satisfies any oracle chain condition. It is therefore easy to use Claim 3.3 on page
127 of \cite{shel.pf} to obtain $\frak{M}_{\xi}$ exactly as in the
proof of Theorem
\ref{not6}.

Let $G$ be $\PP_{\omega_2}$ generic over $V[H]$. 
Exactly as in the proof of Theorem \ref{not6}, it can be shown that the 
cofinality 
of any ultrapower of the integers is $\omega_2$ while $\frak{b}
 =\aleph_1$ in $V[G,H]$. On the other hand, the fact that
$\{c_{\xi} : \xi\in\omega_2\}$ is unbounded
 follows from genericity and the fact that $\SS(\xi + 1) = 
\SS(\xi)\ast \DD(\{c_{\zeta} : \zeta\in\xi\})$ --- 
so $c_{\xi}$ is not dominated by
any function from $V[G_{\xi},H_{\xi}]$ where $G_{\xi}$ 
is the restriction of $G$
to $\PP_{\xi}$.
\stopproof
\begin{corol}
Axiom~6 does\label{c6not5.5} not imply Axiom~5.5.
\end{corol}
To see that the model constructed in Theorem \ref{6not5.5} is a model
of Axiom~6 but not of Axiom~5.5 observe first that Axiom~5.5 fails because
$\frak{b} = \aleph_1$ while
the cofinality of any ultrapower of the integers is $\omega_2$. On
the other hand, there is $(\leq^*,\fomom)$-chain
 of length $\omega_2 \geq \frak{d}$.
Using Lemma \ref{integral} it is possible to construct from this 
a $(<^*,\cal{S})$-chain, $\cal{C}'$, of nondecreasing functions.
From Lemma  \ref{suff6} it follows that
there is an ultrafilter $\cal{U}$ on $\omega$ such that $\cal{C}'$ is
cofinal in the ultrapower of the integers modulo $\cal{U}$.
This is the statement of Axiom~6.\stopproof

The partial order $\SS(\lambda)$ can be modified to yield a model where
Axiom~4 holds yet Axiom~5.5 fails. 
Recall that Lemma~\ref{integral} implies that to 
do this it is only necessary to find a model of Axiom~$4(\leq,\fomom)$
and the failure of Axiom~5.5.
\begin{theor}
    If set theory \label{4not5.5} is consistent 
then there is a model of set theory where
 $\frak{b} = \aleph_1$ yet  for
every ultrafilter $\cal U$ there is a $(\leq^*,\fomom)$-chain
of length $\omega_2$ which is cofinal
in $\fomom/\cal U$.
\end{theor}
\proof It will be shown that, assuming $\lozenge_{\omega_2}(\omega_1)$,
 there is a locally Cohen partial order $\PP$ such that if $G$ is
$\PP$ generic then $\cof(\fomom/\cal{U}) = \omega_2$ for
every ultrafilter $\cal U$ in $V[G]$. The fact that $\PP$ is locally
Cohen will guarantee that $\frak{b} = \aleph_1$ in $V[G]$.

To construct $\PP$ some preliminary bookkeeping is required.
 Let $\{D_{\xi} : \xi\in\omega_2\}$
be a $\lozenge_{\omega_2}(\omega_1)$ sequence and let
$\{g_{\xi} : \xi \in\omega_2\}$ enumerate names for elements of $\fomom$ 
which  arise from countable chain condition
forcing partial orders on $\omega_2$.
Also, if $\QQ$ is any partial order of 
size $\aleph_2$ and satisfying the countable chain
condition then any subset of the reals in a $\QQ$ generic extension
has a name of size $\aleph_2$. Consequently, 
it is possible to use subsets of $\omega_2$ 
to code such names for sets of reals. If $\cal X$ is
some name --- in a suitable partial order ---
for a subset   of $[\omega]^{\aleph_0}$ then $c(\cal{X})$
will denote the subset of $\omega_2$ which codes it while if
$X\subseteq \omega_2$ then $d(X)$ will denote the name it codes. The
details of the coding will not be important. Define a partial order
$\prec$ on $\omega_2$ by $\xi\prec \eta$ if and only 
if $c(D_{\xi}) = c( D_{\eta}\cap\xi)$.

Now construct $\PP$ as a finite support iteration of $\{\PP_{\xi} :
\xi\in\omega_2\}$ such that $\PP_{\xi + 1} = \PP_\xi* \CC_{\xi}\times
\DD_{\xi}$ where
$\CC_{\xi}$ adds a Cohen real, $A_{\xi}:\omega\rightarrow 2$ and
$\DD_\xi$ is some partial order which has yet to be defined. 
At the same time, construct a  partial function
$\Theta : \omega_2\rightarrow \omega_2$  so that 
 if $\mu$ is in the domain of $\Theta$ then
  $1\forces{\PP_{\mu}}{d(D_{\mu})\mbox{ is an ultrafilter}}$
and $\Theta(\mu)$ is the minimum ordinal such that
$\Theta(\mu)\notin\{\Theta(\gamma) : \gamma\prec \mu\}$ and such that
  $g_{\Theta(\mu)}$ is a $\PP_{\mu}$ name.

Given $\PP_{\eta}$, define 
 $\DD_\eta$ by $p\in \DD_{\eta}$ if and only if
\begin{itemize}
\item $p = (f_p,\Gamma_p)$
\item $f_p\in\wfomom$
\item $\Gamma_p\in [\eta]^{<\aleph_0}$
\item if $\gamma \in \Gamma_p$ then $\gamma \prec \eta$
\end{itemize}
and $p \leq q$ is defined to hold if and only if
\begin{itemize}
\item $f_p\subseteq f_q$
\item $\Gamma_p\subseteq \Gamma_q$
\item if $\gamma \in \Gamma_p$,
 $A_{\gamma}^{-1}\{k\}\in d(D_\eta)$,
$A_{\gamma}(n) = k$ and
$n\in \dom(f_q\setminus f_p)$
then $f_q(n) \geq \max\{g_{\Theta(\gamma)}(n),F_\gamma(n)\}$
where, for any $\xi \in \omega_2$, 
$F_{\xi}$ is the  generic function added by the partial order
$\DD_{\xi}$ --- to be precise, $F_\xi = \cup\{f_p : p\in G\}$ where 
$G$ is $\DD_\xi$ generic . 
\end{itemize}
If
it is possible to extend $\Theta$  to include $\eta$ in
its domain then do so  --- there is no
ambiguity here because an extension, if it  exists, is unique. 
 
Let $\PP = \PP_{\omega_2}$. It will soon be shown that $\PP$
satisfies the countable chain condition. However,
first suppose that $G$ is $\PP$ generic over $V$ 
and that $\cal U$ is the $\PP$ name for an ultrafilter
in $V[G]$. There is then a stationary set, $S(\cal{U})$,
such that if $\xi \in S(\cal{U})$ then
 $1\forces{\PP_{\xi}}{d(c(\cal{U})\cap \xi) \mbox{ is an
ultrafilter}}$.
It will be shown that $\{F_{\xi} : \xi\in S(\cal{U})\}$ is a
$(\leq^*,\fomom)$-chain which is cofinal in $\fomom/\cal U$. The
fact that it  is a
$\leq^*$ increasing sequence is an immediate consequence of the
definition of $\DD_\xi$. 

To see that it
 is cofinal in $\fomom/\cal U$ let $g\in\fomom$. Then, assuming that
$\PP$ has the countable chain condition, there is some $\theta\in\
\omega_2$ and $\mu \in\omega_2$
such that $g_{\mu}$ is a $\PP_\theta$ name for $g$. It follows 
 that there is some  $\zeta \in S(\cal{U})$ such that $\mu =
\Theta(\zeta)$. Let $\eta\in S(\cal{U})\setminus(\zeta + 1)$
and note that $\zeta \prec \eta$.
Hence, the partial order $\DD_\eta$ adds a
function which dominates $g_\mu$ on $A_\zeta^{-1}\{k\}$ for some $k\in
2$
and, moreover, $A_\zeta^{-1}\{k\}\in\DD_\eta$.

It remains to be shown that $\PP$ satisfies the countable chain condition
 and that $\frak{b} = \aleph_1$ after forcing with $\PP$. Both these facts
will follow once it has been  shown  that $\PP$ is locally
Cohen. To this end, it is worth observing that $\PP$ has a dense set
of conditions which are somewhat determined  --- a condition $p$ will
be said to be
{\em somewhat determined}
if the support of $p$ is $\Sigma^p\in [\omega_2]^{<\aleph_0}$ and
there is an integer $n(p)$ such that
\begin{itemize}
\item for each $\sigma\in \Sigma^{p}\cap \dom(\Theta)$ there is 
$h_0^{p,\sigma} :n(p)\rightarrow 2$,  $h_1^{p,\sigma}:n(p)\rightarrow
\omega$
 and
$\Delta^{p,\sigma}\in [\sigma]^{<\aleph_0}$ such that
 $p\restriction \sigma\forces{\PP_{\sigma}}{p(\sigma) = (h_0^{p,\sigma},
h_1^{p,\sigma},\Delta^{p,\sigma})} $
\item $\Delta^{p,\sigma} \subseteq  \{\zeta \in\Sigma^{p}\cap\sigma
: \zeta \prec\sigma\}$ 
\item for each $\sigma\in \Sigma^p$ and $\tau \in \Sigma^p$
 such that $\sigma \prec \tau$  there
is
$k(p,\sigma,\tau)\in 2$ such that
$p\restriction \tau\forces{\PP_\tau}{A_{\sigma}^{-1}\{k(p,\sigma,\tau)\}\in 
d(D_{\tau})}$
\item for each $\sigma \in \Sigma^p$ and  $\tau \in \Sigma^p$
such that $\sigma \prec \tau$ there is $M_p(\sigma)\in \omega$ and $G_p(\sigma) :
M_p(\sigma)\rightarrow \omega$ such that $p\restriction
\sigma\forces{\PP_{\sigma}}{g_{\Theta(\sigma)}\restriction M_p(\sigma) = G_p(\sigma)}$
and, moreover, $h_0^{p,\sigma}(i) = 1 + k(p,\sigma,\tau)\mod 2$
 provided that $i\in n(p)
\setminus
 M_p(\sigma)$.
\end{itemize}
 The fact that the set of somewhat determined conditions in $\PP_\eta$
 is dense in $\PP_\eta$ will be
 proved by induction, but an extra induction hypothesis is necessary. 
What will be shown by induction on $\eta$ is that, given 
\begin{itemize}
\item $p\in \PP_\eta$
\item any finite set $W$ of maximal elements of $\prec\cap(\eta\times \eta)$
\item any function $v : W\rightarrow 2$
\item any function $a:W \rightarrow \omega_2$ such that $g_{a(\xi)}$ is
a $\PP_\xi$ name for each $\xi \in W$
\end{itemize}
there is a determined condition $q$ --- the fact that $q$ is
determined is witnessed by $n(q)$ --- with the additional properties
that
for each $\xi \in W$ there is $M(\xi)\in \omega$ and $G(\xi) :
M(\xi)\rightarrow \omega$ such that $q\restriction
\xi\forces{\PP_{\xi}}{g_{a(\xi)}\restriction M(\xi) = G(\xi)}$
and, moreover, $h_0^{p,\xi}(i) = v(\xi)$ provided that $i\in n(q)
\setminus
 M(\xi)$.

If $\eta = 0$  this is trivial and
if $\eta$ is a limit ordinal then  it follows from the fact that a finite
support iteration is being used. 
Therefore, suppose that the fact has been established for $\eta $
 and that $p =
(p\restriction \eta, (h_0, h_1, \Gamma))\in \PP_{\eta + 1}$. Suppose
also that $W$, 
 $v : W\rightarrow 2$ and $a:W \rightarrow \omega_2$ have been given
so that $W$  is a finite set
of maximal elements of $\prec\cap(\eta + 1)^2$.
Notice that $\prec\cap(\eta + 1)^2$ has at most one  maximal
element, $\eta$,
which is not maximal in $\prec\cap(\eta\times \eta)$. It is, of
course, possible that some maximal element in $\prec\cap(\eta\times
\eta)$
is no longer maximal in $\prec\cap(\eta + 1)^2$.
 If there is  such a
new non-maximal element, then denote it by $\theta$; if not, then the following
argument is a bit easier and so it will be assumed that $\theta$ exists.
Find $q \geq p\restriction\eta$ and $H_0:I_0\rightarrow 2$ and
$H_1:I_1\rightarrow \omega$ as well as $k\in 2$ and 
$\Delta\in [\eta]^{\aleph_0}$ such that 
\begin{itemize}
\item $q\forces{\PP_\eta}{h_0 = H_0 \AND h_1 = H_1}$
\item $q\forces{\PP_{\eta}}{g_{a(\eta)}
\restriction I_0 = G }$ for some  $G: I_0\rightarrow
\omega$ (if $\eta\notin W$ this can be ignored)
\item $q\forces{\PP_\eta}{\Gamma = \Delta}$
\item $q\forces{\PP_{\eta}}{A_{\theta}^{-1}\{k\}
\in d(D_\eta)}$
\end{itemize}
That it is possible to arrange for the first two clauses follows from the
fact that $g_{a(\eta)}$ is a $\PP_{\eta}$ name and so any
information about it can be obtained without changing
$h_0$ or $h_1$. To satisfy the last clause, use the fact that $\theta\prec
\eta$, which follows because $\theta$ is no longer maxinal in $\prec\cap(\eta + 1)^2$.

Now define
 $W' = (W\setminus \{\eta\})\cup \{\theta\}$ and observe that $W'$ is
a set of maximal elements in $\prec\cap(\eta\times \eta)$.
Define $v' = v\restriction W'\cup\{(\theta,k+1\mod 2)\}$ and
$a' = a\restriction W'\cup\{(\theta,\Theta(\theta))\}$ and observe that
 both $a'$ and $v'$
are still functions of the right type.
 Then use the induction hypothesis on $\eta$ to find $q' \geq q$ which
is somewhat determined and such that this is witnessed by $n(q')$ and,
such that
for each $\xi \in W'$ there is $M(\xi)\in \omega$ and $G(\xi) :
M(\xi)\rightarrow \omega$ such that $q'\restriction
\xi\forces{\PP_{\xi}}{g_{a'(\xi)}\restriction M(\xi) = G(\xi)}$
and, moreover, $h_0^{p,\xi}(i) = v'(\xi)$ provided that $i\in n(q)
\setminus
 M(\xi)$.
Without loss of generality, $n(q') \geq I_0$ and $n(q') \geq I_1$.

Then let $p' =(q', (h_0',h_1', \Gamma)$ where $h_0':n(q')\rightarrow
2$ is the extension of $H_0$ to $n(q')$ such that $h_0'(i) = v(\eta)$
if $i\in n(q')\setminus I_0$ and $h_1'$ is the extension of
$H_1$ such that 
$$h_1'(i) = \max(\{f_\gamma(i) : \gamma\prec \eta\AND
\gamma \in \Gamma\}\cup\{G(\Theta(\gamma))(i) : \gamma\prec \eta\AND
\gamma \in \Gamma\})$$
 for $i\in n(q')\setminus I_1$. Notice that maximum is taken over 
actual integers rather than names for integers. The definition also
respects the requirements of extension in the partial order
$\PP_\eta$.
Defining $M(\eta) = I_0$ and $G(\eta)= G$ satisfies the extra
induction hypothesis.

To see that $\PP$ is locally Cohen 
let $X\in [\PP]^{\aleph_0}$. Let $\frak{M}$
be a countable elementary submodel of $(H(\omega_3),
\PP,\{D_{\xi} : \xi\in\omega_2\},\Theta,X)$. It suffices to show that
$\PP\cap \frak{M} $ is completely embedded in $\PP$. To see that
this is so, let $A\subseteq \PP\cap \frak{M}$ be a maximal antichain in
$\PP\cap \frak{M}$ and let $p\in \PP$; without loss of generality
$p$ can be assumed to be somewhat determined and, moreover, it may
be assumed that  this
is witnessed
by $n(p)$. Let $p'$ be defined so that $\dom(p') =\dom( p)\cap \frak{M}$
and  $p'(\xi) = (h_0^{p,\xi}, h_1^{p,\xi},\Sigma^{p}\cap\xi\cap
\frak{M})$
 for $\xi
\in\dom(p')$. Note that $p'\in\frak{M}\cap \PP$.
Hence there is $q'\in A$ and $q\in \PP\cap \frak{M}$ such that
$q\geq q'$ and $q \geq p'$ --- without loss of generality it may
be assumed that $q$ is determined and this is witnessed by $n(q)$.
 It must be shown that $p$ and $q$ are
compatible.
 
As in the proof that $\SS(\lambda)$ is locally Cohen, for $\sigma \in
\dom(p) \setminus \dom(q)$ extend
$h_1^{p,\sigma}$ to $ h_1^\sigma$ by defining
 $$h_1^\sigma(m) = \max(\{ h_1^{q,\tau}(m) : \tau\prec \sigma \AND
 \tau\in \dom (q)\}\cup\ldots$$$$\ldots
\{ G_q(\tau)(m) : \tau\prec \sigma \AND
 \tau\in \dom (q)\AND m\in M_q(\tau) \AND A_\tau(m) \neq k(q,\tau,\sigma)\})$$
 for $m \in n(q) \setminus n(p)$, recalling that $G_q(\xi)$, $M_q(\xi)$
and $k(q,\xi,\eta)$ are witnesses to the fact that $q$ is somewhat
determined.
 This will
certainly
assure that
 if $\tau$ and $\rho$ are in the domain of $q$ and
$\tau \prec \sigma \prec \rho$ then
 $h_1^{q,\tau}(m) \leq h_1^{\sigma}(m)\leq
h_1^{q,\rho}(m)$\ \ --- the fact that $ h_1^{\sigma}(m)\leq
h_1^{q,\rho}(m)$ follows from the defintion of the third coordinates
in $p'$. Also,
$g_{\Theta(\tau)}(m) \leq h_1^{\sigma}(m)$ if $A_\tau(m)\neq
k(p,\tau,\sigma)$ will be true if $m \in M(\tau)$ by construction.

Next, if  $\sigma\in
\dom(p)\setminus \frak M$, $\tau \in \dom(p)\cap \frak M$
 and $\sigma \prec \tau$ then define $h_0^{\sigma}(m)
= 1 + k(p,\sigma,\tau)\mod 2$. If this can be done then, if $m \geq n(p)$,
 it is not necessary for $h_1^{q,\tau}(m)$
to be greater than $ g_{\Theta(\sigma)}(m)$.
If there is no $\tau \in
\dom(p)\cap \frak{M}$ such that $\sigma\prec\tau$, do not extend
$h_0^{p,\sigma}$ at all. Notice that in this last case it is still
possible that there is some
$\tau \in \dom (q)$ such that $\sigma \prec \tau$. However, because
it is only necessary  for $h_1^{q,\tau}(m)$
to be greater than $ g_{\Theta(\sigma)}(m)$ in case $\sigma \in
\Delta^{q,\tau}$, this will cause no problems because
$\Delta^{q,\tau}\subseteq \frak M$ if $\tau\in \frak M$.

 What must be checked, though, is that no 
conflict arises as a result of this definition of
$h_0^{\sigma}$. After all, it is conceivable that
$\sigma \prec \tau$ and $\sigma \prec \tau'$
  but $k(q,\sigma,\tau)\neq k(q,\sigma,\tau')$.
To see that this does not happen,
suppose that $\sigma\prec \tau $, $\sigma\prec\tau'$,
 $k(p,\sigma,\tau)\neq k(p,\sigma,\tau')$ and $\{\tau,\tau'\}
\subseteq \frak{M}$. It follows that if $$\rho =
\sup\{\theta : \theta\prec \tau \AND \theta\prec \tau'\AND
\theta\in\dom(\Theta)\}$$
 then $\rho
\in\frak{M}$. Hence $\sigma \in \rho$ and so there is some
$\theta'$ such that $\theta'\prec \tau$, $ \theta'\prec \tau'$ and
$\theta\in\dom (\Theta)$ such that $\sigma \in \theta'$. 
Hence $A_{\sigma}^{-1}\{0\}$ is measured
by the ultrafilter $d(D_{\theta'})$. Since
$d(D_{\theta'})\subseteq d(D_{\tau})$ and
$d(D_{\theta'})\subseteq d(D_{\tau'})$ it follows that
 $ k(p,\sigma,\tau) = k(p,\sigma,\tau')$.
\stopproof
\section{Open Questions}
Table~\ref{tab1} of implications and non-implications 
summarizes the known results about the axioms discussed in this paper.
 The key to understanding
Table~\ref{tab1} is that
\begin{itemize}
\item if there is a ``$\Rightarrow$'' in the entry in the row headed by 
Axiom $R$
and the column headed by Axiom $C$ then Axiom $R$ implies Axiom $C$
\item if there is a ``$\not\Rightarrow$'' in the entry in the row headed by Axiom $R$
and the column headed by Axiom $C$ then Axiom $R$ is known to be consistent with the negation of Axiom $C$
\item if there is a question mark in the entry in the row headed by Axiom $R$
and the column headed by Axiom $C$ then it is not known whether Axiom $R$ implies Axiom $C$
\end{itemize}

\newcommand{\iii}{$\Rightarrow$}
\newcommand{\nii}{$\not\Rightarrow$}

\begin{table}\caption{Table of Implications}
\begin{center}\label{tab1}
    \begin{tabular}{|c|c|c|c|c|c|c|c|c|c|c|} 
\hline
Axiom          & 1  & 2  & 3  & 4  & 5$\leq^*\fomom$& 5$<^*\fomom$& 5$\cal{N}$& 5.5 & 6  & 6.5  \\ \hline
1              &\iii&\iii&\iii&\iii&     ?          &   ?         &\iii      &\iii  &\iii&\iii  \\ \hline
2              &\nii&\iii&\nii&\iii&  \nii          &   \nii      & \iii     &\iii  &\iii&\iii  \\ \hline
3              &\nii&\nii&\iii&\nii&   \nii         &    \nii     & \iii     &\iii  &\iii&\iii  \\ \hline
4              &\nii&\nii&\nii&\iii&  \nii          &   \nii      & \iii     &\nii  &\iii&\iii  \\ \hline
5$\leq^*\fomom$&\nii& ?  & ?  & ?  &  \iii          &   \iii      & \iii     &\iii  &\iii&\iii  \\ \hline
5$<^*\fomom$   &\nii& ?  & ?  & ?  &  ?             &    \iii     & \iii     &\iii  &\iii&\iii  \\ \hline
5$\cal N$      &\nii&\nii&\nii&\nii& \nii           &   \nii      & \iii     &\iii  &\iii&\iii  \\ \hline
5.5            &\nii&\nii&\nii&\nii&  \nii          &    \nii     & ?        &\iii  &\iii&\iii  \\ \hline
6              &\nii&\nii&\nii&\nii&  \nii          &     \nii    & \nii     &\nii  &\iii&\iii  \\ \hline
6.5            &\nii&\nii&\nii&\nii&  \nii          &     \nii    & \nii     &\nii  & ?  &\iii  \\ \hline
     \end{tabular}\end{center}
\end{table}

Not all the reasons for the assertions made in Table~\ref{tab1} are contained
in in this paper. Some will be found in in \cite{nyik.c} and others must be deduced
by modus ponens. Table~\ref{tab2}
  contains a guide to reasons for the various
 assertions in  Table~\ref{tab1}. 

\begin{table}\caption{Table of References}\label{tab2}
\begin{center}
    \begin{tabular}{|c|c|c|c|c|c|c|c|c|c|c|} 
\hline
Axiom          & 1  & 2  & 3  & 4  & 5$\leq^*\fomom$& 5$<^*\fomom$& 5$\cal{N}$& 5.5 & 6  & 6.5  \\ \hline
1              & T  & T  & T  & T  &  ?             &   ?         &    T      &  N  & T  & N    \\ \hline
2              & N  & T  & 1  & T  &  1             &   1         &  2        &  N  & T  & N    \\ \hline
3              & N  & N  & T  & N  &  9             &   9         &  T        &  N  & T  & N    \\ \hline
4              & N  & N  & N  & T  &  3             &   3        &  T        &  4  & T  & N    \\ \hline
5$\leq^*\fomom$& 10 & ?  & ?  & ?  &  T             &   T         &  T        &  N  & T  & N    \\ \hline
5$<^*\fomom$   & 10 & ?  & ?  & ?  &  ?             &   T         &  T        &  N  & T  & N    \\ \hline
5$\cal N$      & N  & N  & N  & N  &  9             &   9         &  T        &  N  & T  & N    \\ \hline
5.5            & 5  & 5  & 5  & 5  &  9             &   9         &  ?        &  T  & N  & N    \\ \hline
6              & 5  & 5  & 5  & 5  &  6             &   6         &  8        &  7  & T  & N    \\ \hline
6.5            & 5  & 5  & 5  & 5  &  6             &   6         &  5        &  5  & ?  & T     \\ \hline
     \end{tabular}\end{center}
\end{table}

A ``T'' in the row  corresponding
to Axiom~$R$ and the column corresponding to Axiom~$C$ in Table~\ref{tab2}
indicates that the
fact that Axiom~$R$ implies Axiom~$C$ is a trivial implication --- 
by trivial is meant
 something which can be deduced by considering the quantifiers in the 
relevant axioms. An ``N'' in that
entry means that either the implication or non-implication can be found in
\cite{nyik.c}. 
The enumeration of the following list corresponds to the numbered entries in 
Table 2. So, for example, the second entry of this list refers 
to  Theorem~\ref{2implies5}
because this is the reason there is an ``\iii''
 in Table 1 in the row corresponding
to Axiom~2 and the column corresponding to Axiom~$5(\cal{N})$.
\begin{enumerate}
\item Theorem~\ref{CHth}
\item Theorem~\ref{2implies5} 
\item The fact that Axiom~4
does not imply Axiom~$5(<^*,\fomom)$ follows because
it has been shown that Axiom~2 does not imply Axiom~$5(<^*,\fomom)$
in Theorem~\ref{CHth} and the fact that Axiom~2 implies Axiom~4 follows
from an inspection of the quantifiers involved. Modus ponens
yields the rest. 
\item Theorem~\ref{4not5.5}
\item The antecedent of the implication is implied by Axiom~$5(\cal{N})$ so 
the non-implication follows from modus ponens because Axiom~$5(\cal{N})$ 
does not imply the conclusion.
\item The antecedent of the implication is implied by Axiom~$4$ so 
the non-implication follows from modus ponens.
\item Corollary~\ref{c6not5.5}
\item Axiom~$5(\cal{N})$ implies Axiom~5.5.
\item $2^{\aleph_0}= \aleph_1$ is known \cite{dow.rempt}
 to imply Axiom~3 and Lemma \ref{CHlem} shows that
Axiom~$5(<^*,\fomom)$ fails under this assumption. 
For the rest, use modus ponens.
\item See the remarks following Theorem~\ref{pfaimp5}. 
\end{enumerate}
%\bibliography{c:/os2/apps2/emtex/bibinput/myabbrev,c:/os2/apps2/emtex/bibinput/standard}
\bibliography{myabbrev,standard}
\end{document}